\documentclass[a4paper,english]{extarticle}
\usepackage[T1]{fontenc}
\usepackage[utf8]{luainputenc}
\usepackage{float}
\usepackage{amssymb}

\makeatletter

\pdfpageheight\paperheight
\pdfpagewidth\paperwidth

\providecommand{\tabularnewline}{\\}

\newcommand{\lyxrightaddress}[1]{
\par {\raggedleft \begin{tabular}{l}\ignorespaces
#1
\end{tabular}
\vspace{1.4em}
\par}
}

\makeatother

\usepackage{babel}
\begin{document}

\title{Primes of the form {\Large{}$p=1+n!\sum n,$ for some $n\in\mathbb{N}^{+}$}}

\author{{\normalsize{}By}\\
Maheswara Rao Valluri{\normalsize{} }}
\maketitle
\begin{abstract}
The purpose of this note is to report on the discovery of the primes
of the form $p=1+n!\sum n$, for some natural numbers $n>0$. The
number of digits in the prime p are approximately equal to $\lfloor log_{10}(1+n!\sum n)\rceil+1$.
\\


\end{abstract}

\section*{Introduction}

A natural number is a prime if it has only factors of $1$ and itself.
There are, by Euclidean theorem (about 350BC) infinitely many primes.
There are many patterns of primes amongst which the classical known
primes are the Mersenne primes of the form $2^{p}-1$, where $p$
is a prime {[}Da2011{]}, and the Fermat primes of the form $2^{2^{n}}+1$
for a natural number $n\geq0$ {[}Da2011{]}. We omit the other classes
of primes except the factorial primes. 

There was a known fact that early 12th century Indian scholars knew
about factorials. In 1677, a British mathematician, Fabian Stedman
described factorials for music. In 1808, a French mathematician, Christian
Kramp introduced notation $!$ for factorials. The factorial of $n$
can be described as product of all positive integers less than or
equal to $n$. In the Christian Kramp's notation, $n!=n(n-1)(n-2)......3.2.1.$
Factorials of $0$ and $1$ can be written as $0!=1$ and $1!=1$,
respectively. There are dozens of prime factorials. However, we recall
only few of them, particularly, factorial primes of the form $(p!\pm1)$
{[}Bor1972, BCP1982{]}, double factorial primes of the form $n!!\pm1$
for some natural number $n$ {[}Mes1948{]}, Wilson primes: $p$ for
which $p^{2}$ divides $(p-1)!+1$ {[}Bee1920{]}, and primorial primes
of the form $(p\#\pm1)$ which means as the product of all primes
up to and including the prime {[}Dub1987, Dub1989{]}. Further, a class
of the Smarandache prime is of the form $n!\times S_{n}(n)+1$ , where
$S_{n}(n)$ is smarandache consecutive sequence {[}Ear2005{]}. 

The purpose of this note is to report on the discovery of the primes
of the form $p=1+n!\sum n$ for some natural numbers $n>0$.

\section*{Primes of the form $p=1+n!\sum n$ for some $n\in\mathbb{N}^{+}$}

We list in the table 1, the primes of the form $p=1+n!\sum n$ for
some natural numbers $n>0$. They are verified up-to $n=950$ and
primes are found when $n=1,2,3,4,5,6,$ $7,8,9,10,12,13,14,19,24$,251,374.
The above primes can also be expressed as $p=1+\frac{(n+1)!n}{2},$
for some natural numbers $n>0$. The author has used python software
to search and verify the above form primes and could verify up-to
$n=950$. The author conjectures that there are infinitely many primes
of the form $p=1+n!\sum n$ for some $n\in\mathbb{N}^{+}$ .

\begin{table}[H]
\caption{Primes List}

~

\centering{}%
\begin{tabular}{|c|c|c|c|}
\hline 
S.No. & n & Prime, $p$ & Number of \tabularnewline
 &  &  & digits in $p$\tabularnewline
\hline 
\hline 
1 & 1 & 2 & 1\tabularnewline
\hline 
2 & 2 & 7 & 1\tabularnewline
\hline 
3 & 3 & 37 & 2\tabularnewline
\hline 
4 & 4 & 241 & 3\tabularnewline
\hline 
5 & 5 & 1801 & 4\tabularnewline
\hline 
6 & 6 & 15121 & 5\tabularnewline
\hline 
7 & 7 & 141151 & 6\tabularnewline
\hline 
8 & 8 & 1451521 & 7\tabularnewline
\hline 
9 & 9 & 16329601 & 8\tabularnewline
\hline 
10 & 10 & 199584001 & 9\tabularnewline
\hline 
11 & 12 & 37362124801 & 11\tabularnewline
\hline 
12 & 13 & 566658892801 & 12\tabularnewline
\hline 
13 & 14 & 9153720576001 & 13\tabularnewline
\hline 
14 & 19 & 23112569077678080001 & 20\tabularnewline
\hline 
15 & 24 & 186134520519971831808000001 & 27\tabularnewline
\hline 
16 & 251 & 25662820338985371726..Omitted..000000000001 & 500\tabularnewline
\hline 
17 & 374 & 22873802587990440054..Omitted..0000000000001 & 807\tabularnewline
\hline 
\end{tabular}
\end{table}

\subsection*{Size of the prime of the form $p=1+n!\sum n$ for some $n\in\mathbb{N}^{+}$}

To compute the size of primes of above form, we use the Stirling 's
formula {[}Sec 2.2, CG2001{]}: $log\,n!=(n+\frac{1}{2})log\,n-n+\frac{1}{2}log2\pi+O(\frac{1}{n})$.
Simply, we can also write $log\,n!\sim n(logn-1)$, if necessary.
The size of the prime, $p=1+n!\sum n$ is approximately equal to $\lfloor log_{10}(1+n!\sum n)\rceil+1$.

\section*{Conclusion }

In this note, the author conjectures that there are infinitely many
primes of the form $p=1+n!\sum n$ for some natural numbers $n>0$.
The author has found the primes for $n=1,2,3,4,5,5,6,7,8,9,10,12,13,14,19,24,251,374$,
when they are searched and verified up-to $n=950$. The number of
digits in the primes of the form $(1+n!\sum n)$ are approximately
equal to $\lfloor log_{10}(1+n!\sum n)\rceil+1$. Furthermore, an
investigation will be required for finding such primes.

\lyxrightaddress{School of Mathematical and Computing Sciences\\
Fiji National University\\
P.O.Box:7222, Derrick Campus, Suva, Fiji\\
E-mail: maheswara.valluri@fnu.ac.fj}
\end{document}